\documentclass[notitlepage,reqno,11pt]{amsart}
\usepackage{latexsym,amssymb, epsfig, amsmath,amsfonts, subfigure,amsthm, mathrsfs}

\usepackage{rotating}
\usepackage[toc,page]{appendix}
\usepackage{color}
\usepackage{multirow}
\usepackage{relsize}
\usepackage{microtype}
\usepackage{dsfont}
\usepackage{cases}
\usepackage[foot]{amsaddr}
\usepackage[colorlinks,anchorcolor=blue,citecolor=blue,allcolors=blue]{hyperref}
\usepackage{enumitem}
\usepackage[colorlinks, allcolors=blue]{hyperref}
\usepackage{xpatch}  
\makeatletter   
\xpatchcmd{\@thm}{\thm@headpunct{.}}{\thm@headpunct{}}{}{}   

\usepackage[margin=1in]{geometry}

\usepackage{xpatch}  
\makeatletter   
\xpatchcmd{\@thm}{\thm@headpunct{.}}{\thm@headpunct{}}{}{}

\numberwithin{equation}{section}
\newtheorem{theorem}{Theorem}[section]
\newtheorem{lemma}{Lemma}[section]
\newtheorem{assumption}{Assumption}[section]

\newtheorem{proposition}{Proposition}[section]
\newtheorem{remark}{Remark}[section]

\newenvironment{pr}{{\bf{Proof.~}}}

\newlength{\defbaselineskip}
\newcommand{\setlinespacing}[1]%
{\setlength{\baselineskip}{#1 \defbaselineskip}}

\newcommand{\Z}{{\mathbb Z}}

\def\E{\mathbb{E}}
\def\P{\mathbb{P}}

\newcommand{\eps}{\varepsilon}

\newcommand{\R}  {\mathbb{R}}
\newcommand{\N}  {\mathbb{N}}

\newcommand{\baa}{\begin{eqnarray*}}
	\newcommand{\eaa}{\end{eqnarray*}}

\def\eps{\varepsilon}

\def\P{{\mathbb P}}
\def\E{{\mathbb E}}
\def\Z{{\mathbb Z}}
\def\N{{\mathbb N}}
\def\R{{\mathbb R}}

\def\1{{\mathbf 1}}

\def\eps{\varepsilon}

\newcommand{\tor}{\mathbb{T}^d}

\newcommand{\xp}{x_\eps}
\newcommand{\yp}{y_\eps}
\newcommand{\bnep}{B^{N, \eps}}
\newcommand{\snep}{S^{N, \eps}}
\newcommand{\inep}{I^{N, \eps}}
\newcommand{\fnep}{\mathfrak{F}^{ N, \eps}}

\newcommand{\ovsnp}{\overline{S}^{\, N, \eps}}
\newcommand{\ovinp}{\overline{I}^{\, N, \eps}}
\newcommand{\ovfnp}{\overline{\mathfrak{F}}^{\, N, \eps}}

\newcommand{\dis}{\displaystyle}
 \newcommand{\vx}{V_{\eps}(\xp)}

 \newcommand{\mds}{\mathds{1}}
 \newcommand{\n}{\nonumber}
 
 \newcommand{\sep}{\overline{S}^{\,\eps}}
 \newcommand{\iep}{\overline{I}^{\,\eps}}
 \newcommand{\bep}{\overline{B}^{\,\eps}}
 \newcommand{\xep}{\overline{X}^{\,\eps}}
 \newcommand{\fep}{\overline{\mathfrak{F}}^{\,  \eps}}

 \newcommand{\sepn}{\overline{S}^{\, N, \eps}}

 \newcommand{\xepn}{\overline{X}^{\, N, \eps}}
 \newcommand{\fepn}{\overline{\mathfrak{F}}^{\, N, \eps}}
 
 \newcommand{\sbf}{\overline{\mathbf{S}}^{\,\eps}}
 \newcommand{\ibf}{\overline{\mathbf{I}}^{\,\eps}}
 
 \newcommand{\xbf}{\overline{\mathbf{X}}^{\,\eps}}
 \newcommand{\fbf}{\overline{\mathbf{F}}^{\,\eps}}
 
 \newcommand{\ssbf}{\overline{\mathbf{S}}^{\, N, \eps}}
 \newcommand{\sibf}{\overline{\mathbf{I}}^{\, N, \eps}}
 \newcommand{\sbbf}{\overline{\mathbf{B}}^{\, N, \eps}}
 \newcommand{\sxbf}{\overline{\mathbf{X}}^{\, N, \eps}}
 \newcommand{\sfbf}{\overline{\mathbf{F}}^{\, N, \eps}}
 
 \newcommand{\sbff}{\overline{\mathbf{S}}}
 \newcommand{\ibff}{\overline{\mathbf{I}}}
 \newcommand{\fbff}{\overline{\mathbf{F}}}
 \newcommand{\bbff}{\overline{\mathbf{B}}}
 \newcommand{\xbff}{\overline{\mathbf{X}}}
 
\newcommand{\wi}{\widetilde{\mathscr{M}}^{\;N,\eps}_{I}}
\newcommand{\ws}{\widetilde{\mathscr{M}}^{\;N,\eps}_{S}}
\newcommand{\msi}{\mathscr{M}^{N,\eps}_{SI}}
 
 \newcommand{\fpr}{\begin{flushright}
 		$\square$
 \end{flushright}}

\newcommand{\ttl}{\Large 
Epidemic models with varying infectivity on a refining spatial grid.\;
I. The SI model}

\begin{document}
	
\title[]{\ttl}
	
\author[Anicet Mougabe-Peurkor]{Anicet Mougabe-Peurkor$^1$}
\address{$^1$Laboratoire de Mathématiques et applications, Université Félix Houphouët Boigny, Abidjan, Côte d'Ivoire}
\email{mougabeanicet@yahoo.fr}

\author[{\'E}tienne Pardoux]{{\'E}tienne Pardoux$^2$}
\address{$^2$Aix Marseille Univ, CNRS, I2M, Marseille, France}
\email{etienne.pardoux@univ.amu.fr}

\author[T\'enan Yeo]{T\'enan Yeo$^1$}
\email{yeo.tenan@yahoo.fr}

\date{\today}

\maketitle

\allowdisplaybreaks

\newcommand{\dispace}{\mathrm{D}_\eps}

\begin{abstract}
We consider a space-time SI  epidemic model with infection age dependent infectivity and non-local infections constructed on a grid of the torus $\tor =(0, 1]^d$, where the individuals may migrate from node to another. The migration processes in either of the two states are assumed to be Markovian. We establish a functional law of large numbers by letting jointly $N$ the initial approximate number of individuals on each node go to infinity and  $\eps$ the mesh size of the grid go to zero. The limit is a system of parabolic PDE/integral equations. The constraint on the speed of convergence of  the parameters $N$ and $\eps$ is that $N \eps^d \to \infty$ as $(N, \eps)\to (+\infty, 0)$. 
\end{abstract}


\section{Introduction}

We consider an epidemic model on a refining grid of the $d$ dimensional torus $\mathbb{T}^d$. Like in the earlier work 
\cite{MTE20}, the individuals move from one patch to its neighbors according to a random walk. The first novelty of this paper is that the infectivity of each individual is a random function, which evolves with the time elapsed since infection, as first considered in \cite{KM27}, and recently studied in \cite{FPP21} and \cite{FPP22}. The second novelty is that we allow infection of a susceptible individual by infectious individuals located in distinct patches, and we use a very general rate of infections. 

There are two  parameters in our model, $N$ which is the order of the number of individuals in each patch, and $\eps$, which is the distance between two neighboring sites. The total number of patches is $\eps^{-d}$, and the total number of individuals in the model is $N\eps^{-d}$. Our goal is to study the limit of the renormalized stochastic finite population model as both $N\to\infty$ and $\eps\to0$. In this paper we obtain a convergence result in $L^\infty$ under the restriction that $N\eps^d\to\infty$.  In \cite{MTE20}, the restriction was much weaker, thanks to clever martingale estimates due to Blount \cite{B}. However, in contradiction with the model in \cite{MTE20}, our model is non Markovian, and several of the fluctuating processes are not martingales. As a result, it does not seem possible to extend the techniques of \cite{B} to the situation studied in the present paper. 

There are three models in the present paper. The stochastic SDE model parametrized by the pair $(N,\eps)$, the deterministic model which is an ODE parametrized by $\eps$ on the patches (and is the LLN limit of the first model when $N\to\infty$ with $\eps$ fixed), and the PDE model on the torus $\mathbb{T}^d$, which is the limit of the ODE model as $\eps\to0$.
The convergence of the ODE model to the PDE model exploits standard arguments on semigroup  and their approximation, based on some result in \cite{TKTO}. The main new argument in the present paper consists in showing 
that the difference in $L^\infty$ between the stochastic and the ODE models, which tends to zero as $N\to\infty$ while $\eps$ is fixed according to \cite{FPP22}, tends also to zero when $(N,\eps)\to(+\infty,0)$, provided $N\eps^d\to \infty$.

In this paper, we consider the SI model, S as susceptible, I as infected. An infected individual has an age of infection dependent infectivity, which we suppose to vanish after some random time. It would be natural to decide that at that time the individual leaves the I compartment, and enters the R compartment, R as recovered. For the sake of simplifying our model, we decide that after being infected, an individual remains in the I compartment for ever. This does not affect the evolution of the epidemic, since when its infectivity remains zero, an individual does not contribute anymore to the propagation of the illness, exactly as an individual in the R compartment of an SIR model. However, there are two drawbacks of the present model. First, we do not follow the evolution of the number of infectious individuals, since we have so to speak merged the I and the R compartments. Second, while we distinguish the rate of movements of the S type and the I type individuals we do not distinguish that rate between the infectious and the recovered individuals.
The reason for studying the SI model separately is that, in our ``Varying Infectivity'' model the techniques for proving the convergence as $\eps\to0$ of the ODE model to the PDE model which we are using in the SI case will not be available in the SIR case. One is forced to use different techniques. We will study the extension of the present results to the SIR model in a future publication. But our conviction is that it is worth to present the results in the SI case, due to the possibility in this case of using classical semigroup techniques.

Let us finally comment on the assumptions on the age of infection dependent infectivity. We assume that to each individual who gets infected is attached a random infectivity function, the functions attached to the various individuals being i.i.d., all having the law of a random function $\lambda$ (the law is different for the initially infected individuals). In this paper, as in \cite{FPP22}, we only assume that $\lambda$ belongs a.s. to the Skorohod space of c\`al\`ag function $\mathbf{D}$,  and satisfies $0\le \lambda(t)\le\lambda^\ast$, for some $\lambda^\ast>0$. This is weaker than the assumptions made in \cite{FPP21}. The proof in \cite{FPP22} is quite different from the proof in \cite{FPP21}. Here we use a proof similar to that in \cite{FPP21}. The limitation is that we obtain only the pointwise convergence of the renormalised total infectivity function, while we obtain uniform in $t$ convergence of the proportions of susceptible and infected individuals. We believe that this proof is interesting, due to its simplicity.  
\smallskip

Note that there is some literature on similar models, but mainly without movements of the various individuals, see in particular \cite{AP} for a SIS Markov model, and \cite{PP} for a SIR varying infectivity model. Our previous publication \cite{MTE20} treats a Markov SIR model with movements and only local infections. 

\smallskip

The paper is organized as follows. We describe our model in detail in section 2, in particular the complex form of the rate of infection. In section 3, we state the law of large numbers limit as $N\to\infty$, with $\eps$ fixed, referring to \cite{FPP22} for the proof. In section 4, we take the limit as $\eps\to0$ in the ODE model. Finally, in section 5, 
we study the difference between the stochastic and the ODE model, as $(N,\eps)\to(+\infty, 0)$, and conclude our main result.

\section{Model description}
We consider a total population size  $N\eps^{-d}$ initially distributed on the $\eps^{-d}$ nodes of a refining  spatial grid $\dispace :=[0,1)^d\cap\eps \Z^d$, in which an infection is introduced. Here $\eps$ is the mesh size of the grid (we assume that $\eps^{-1}\in \N\backslash\{0\}$). We focus our attention to the periodic boundary conditions on the hypercube $ [0, 1]^d$, that is, our domain is the torus $ \tor :=[0, 1]^d$. Our results can be extended to a bounded domain of $\R^d$ with smooth boundary, and Neumann boundary conditions.

\bigskip

\subsection{Set-up and notations}~\\

 We split the population in two subsets $S^{N,\eps}$ and $I^{N,\eps}$. $S^{N,\eps}$ stands for the susceptible individuals, who do not have the disease and who can get infected, while  $I^{N,\eps}$ is referred  to the subset of those individuals who are  suffering from the illness  or have recovered from the disease. \\
We shall denote by $\xp$ the nodes of the grid $\dispace$. $S^{N,\eps}(t, \xp)$ denotes the number of susceptible individuals at site $\xp$ at time $t$. Let $\bnep(t, \xp)$ be the total number of individuals at site $\xp$ at time $t$, i.e. $\bnep(t, \xp):= \snep(t, \xp)+\inep(t, \xp)$. We define $\snep(t)$ (resp. $\inep(t)$) as the total number of susceptible individuals (resp. infected individuals) at time $t$ in the whole population, that is: 
$$\snep(t):=\sum_{\xp}\snep(t, \xp), \quad \text{and}\; \;  \inep(t):=\sum_{\xp}\inep(t, \xp) \, , \; \;  \forall t\ge 0.$$
We have $\dis  \bnep(t):=  \sum_{\xp}\bnep(t)=N\eps^{-d} \, , \; \;  \forall t\ge 0.$\\

To each individual $j$ is attached a random infection-age dependent  infectivity process $\{\lambda_{-j}(t)\,  : t\ge  0\}$ or $\{\lambda_j(t)\,  : t\ge  0\}$.   $\lambda_{-j}(t)$  is the infectivity at time $t$ of the $j$-th initially infected individual. The initially susceptible individual $j$ who is infected at a  random time $\tau_j^{N, \eps}$,  has at time $t$ the infectivity $\lambda_j(t-\tau_j^{N,\eps})$, i.e. $\lambda_j(t)$ is the infectivity at time $t$ after its time of infection of the $j$-th initially susceptible individual. 
We assume that $\lambda_j=0$ on $\R_{-}$\, and that $\{\lambda_{-j}\, :  j\ge 1\}$ and $\{\lambda_j\,  : j\ge 1\}$ are  two mutually independent sequences of i.i.d $\R_+$-valued random functions.

\medskip

We define the infected periods of newly and initially infected individual $j>0$ and $j<0$ respectively,  by the random variables $$\eta_{j}:=\sup\{t>0 : \; \lambda_j(t) > 0\}, \; j\in \mathbb{Z}\backslash \{0\}.$$
We define $F(t):= \P\left(\eta_1 \le t \right)$, , $F_0(t):=\P\left(\eta_{-1} \le t \right)$, the distributions functions of $\lambda_j$ for $j\ge 1$ and  for $j\le-1$ respectivelly. Let $F^c(t):= 1-F(t)$ and $F^c_0(t) := 1-F_0(t).$ We moreover define $$\overline{\lambda}(t):=\E\left[ \lambda_1(t) \right]\; \; \text{and}\; \overline{\lambda}_0(t):=\E\left[ \lambda_{-1}(t) \right].$$
Note that, under the i.i.d. assumption of the random variables  $\{\lambda_j(.)\}_{j\ge 1}$, the sequence of random variables $\{\eta_j\}_{j\ge 1}$ is i.i.d. Also, the sequence of random variables $\{\eta_j\}_{j \le -1}$ is i.i.d.
\medskip

We assume that susceptible individuals move from patch to patch according to a time-homogenous Markov process $X(t)$ with jump rates $\nu_S/\eps^2$ and transition  function $$ p^{\xp, \yp}_{\eps}(s, t)= \mathbb{P}\left(X(t)=\yp | X(s)=\xp \right),$$ and while infectious individuals move from patch to patch according to a time-homogeneous Markov process $Y(t)$ with jump rates $\nu_I/\eps^2$ and transition  function 
$$q^{\xp, \yp}_{\eps}(s, t)= \mathbb{P}\left(Y(t)=\yp | Y(s)=\xp \right).$$ 
$\nu_S$ and $\nu_I$ are positive diffusion coefficients for the susceptible and infected
 subpopulations, respectively.
We assume that  those movements of the various individuals are mutually independent. 

In addition, we use $X_j^{s,\xp}(t)$ (resp. $Y_j^{s,\xp}(t)$) to denote the position at time $t$ of the individual $j$ if it is susceptible (resp. infected) during the time interval $(s,t)$, and was in location/node $\xp$ at time $s$.

 For all $\xp \in \dispace$, let $\vx$ be the cube centered at the site $\xp$ with volume $\varepsilon^d$. Let $\mathrm{H}^{\varepsilon} \subset L^2(\tor)$ denote the space of real valued step functions that are constant on each cell $\vx$.  
  
$\Delta_{\varepsilon}$ is the discrete Laplace operator defined as follows
$$\Delta_{\varepsilon}f(\xp) = \sum_{i=1}^{d} \varepsilon^{-d}\big[f(\xp+\varepsilon e_i)-2f(\xp)+f(\xp-\varepsilon e_i) \big], \; \;  f\in \mathrm{H}^{\varepsilon}$$ and we define the operators
$\Delta_\eps^S f:= \nu_S\Delta_\eps f$ and $\Delta_\eps^I f := \nu_I\Delta_\eps f$, $f\in \mathrm{H}^{\varepsilon}$. 

\smallskip

$\Delta$ denotes the d-dimensional Laplace operator. Let $T_{S,\eps}$ (resp. $T_{I,\eps}$) be the semigroup acting on $\mathrm{H}^{\eps}$ generated by $\nu_S\Delta_{\eps}$ (resp.  $\nu_I\Delta_{\eps}$). Similary, we denote by  $T_{S}$ (resp. $T_{I}$)  the semigroup acting on $L^2(\tor)$ generated by $\nu_S\Delta$ (resp.  $\nu_I\Delta$).

\subsection{Model formulation}~\\
All random variables and processes are defined on a common complete probability space $(\Omega , \mathcal{F} , \P)$. We consider a SI epidemic model where each infectious individual has an infectivity that is randomly varying with the time elapsed since infection.
We assume that a susceptible individual in patch $\xp$ has contacts with infectious individuals of patch $\yp$ at rate  $\beta_{\eps}^{\xp, \yp}(t)$ at time $t$.

Given a site $\xp$, the total force of infection at each time $t$ at site $x_\eps$ is the aggregate infectivity of all the individuals that are currently infectious in this site:
\begin{align*}
\mathfrak{F}^{N,\eps}(t,x_\eps)&= \sum_{j=1}^{I^{N,\eps}(0)}\lambda_{-j}(t)\mds_{Y_j(t)=x_\eps}\\
&+
\sum_{y_\eps}\int_0^t\int_0^\infty \int_\mathbf{D}\!\int_\mathbf{D} \lambda(t-s)\mds_{u\le S^{N,\eps}(s^-,y_\eps)\overline{\Gamma}^{N,\eps}(s^-,y_\eps)} \mds_{Y^{s,y_\eps}(t)=x_\eps}Q^{y_\eps}(ds, du, d\lambda, dY),
\end{align*}
where $$\overline{\Gamma}^{N,\eps}(t,y_\eps) :=\frac{1}{N^{1-\gamma}[B^{N,\eps}(t,y_\eps)]^\gamma}\sum_{x_\eps}\beta_\eps^{y_\eps,x_\eps}(t)\mathfrak{F}^{N,\eps}(t,x_\eps)$$
is the force of infection exerted on each susceptible individual in patch $x_{\varepsilon}$, and $\left\{Q^{y_{\varepsilon}}, y_{\varepsilon} \in D_{\varepsilon}\right\}$  are i.i.d.  standard Poisson random measures (PRM) on $\R_+^2\times \mathbf{D}^2$ with intensity  $ds\otimes du\otimes d\P_{\lambda}\otimes d\P_Y$. $\mathbf{D}$ denotes the space of càdlàg paths from $\R_+$ into $\R_+$, which we equip with the Skorohod topology. 
 We assume that  $\gamma \in [0, 1]$.
By an abuse of notation, we denote  by $Q^{x_\eps}(ds, du)$ the projection of $Q^{x_\eps}(ds, du, d\lambda, dY)$ on the first two coordinates.
Let, with $\Upsilon^{N,\eps}(t,x_\eps):=S^{N,\eps}(t,x_\eps)\overline{\Gamma}^{N,\eps}(t,x_\eps)$,
\[ A^{N,\eps}(t,\xp) :=\int_{0}^{t}\int_{0}^{\infty}\mds_{u\leq \Upsilon^{N,\eps}(s^-,x_\eps)}Q^{\xp}(ds,du).
\]
 In what follows, $\xp \sim \yp$ means that the nodes $\xp$ and $\yp$ are neighbors (each  point of $\dispace$ has $2d$ neighbors).

\medskip

The epidemic dynamic of the model can be described by the following equations

\begin{equation}\label{ssto}
\begin{aligned}
S^{N,\eps}(t,x_\eps) &= S^{N,\eps}(0,x_\eps)-A^{N,\eps}(t,\xp)-\sum_{y_\eps\sim x_\eps}^{}P_{S}^{x_\eps,y_\eps}\left(\int_{0}^{t}\frac{\nu_{S}}{\eps^2}S^{N,\eps}(s,x_\eps)ds\right)
+\sum_{y_\eps\sim x_\eps}P_{S}^{y_\eps,x_\eps}\left(\int_{0}^{t}\frac{\nu_{S}}{\eps^2}S^{N,\eps}(s,y_\eps)ds\right)\\
I^{N,\eps}(t,x_\eps) &= I^{N,\eps}(0,x_\eps) + A^{N,\eps}(t,x_\eps)
-\sum_{\yp\sim \xp}P_{I}^{\xp,\yp}\left(\int_{0}^{t}\frac{\nu_{I}}{\eps^2}I^{N,\eps}(s,\xp)ds\right)
+\sum_{\yp\sim \xp}P_{I}^{\yp,\xp}\left(\int_{0}^{t}\frac{\nu_I}{\eps^2}I^{N,\eps}(s,\yp)ds\right),
\end{aligned} 
\end{equation}

where $P_{S}^{x_\eps,y_\eps}$, $P_{I}^{x_\eps,y_\eps}$, $\xp\, ,  \yp \in \dispace$  are mutually independent standard Poisson processes.\\ In the above equations
 $P_{S}^{x_\eps,y_\eps}$ (resp.  $P_{I}^{x_\eps,y_\eps}$) is the counting process of susceptible (resp. infected) individuals that migrate from the patch $x_\eps$ to $y_\eps$.
 
 \medskip
 
  In the sequel of this paper we may use the same notation for different constants (we use the generic notations $c$, $C$ for  positive constants). These constants can depend upon some  parameters of the model, as long as these are independent of  $ \varepsilon $ and $ N$,  and we will not necessarily mention this dependence explicitly.  The exact value may change from line to line.

\section{Law of large numbers as $N \to \infty$, $\eps$ being fixed}

We consider  the renormalized model by dividing the number of individuals in each compartment and at each patch by $N$. Hence, we define
$$\ovsnp(t,\xp):=\dfrac{1}{N}\snep(t, \xp), \quad \ovinp(t,\xp) :=\dfrac{1}{N}\inep(t,\xp), \; \text{and} \; \;  \ovfnp(t,\xp) :=\dfrac{1}{N}\fnep(t,\xp).$$

\begin{assumption}\label{ass1}
We make the following assumptions on the initial conditions. We assume that:
\begin{enumerate}[label=\bf (\roman*)]
\item  there exists a collection of positive numbers $\{\,\overline{S}^{\eps}(0, \xp), \; \overline{I}^{\eps}(0, \xp), \; \xp\in \dispace \; \}$ such that $$  \sum_{\xp}\Big[\overline{S}^{\eps}(0, \xp)+\overline{I}^{\eps}(0, \xp)\Big]=\eps^{-d} \; ,$$ and 
$\dis \Big|S^{N,\eps}(0)-N\overline{S}^\eps(0)\Big| \le 1 \; , \quad \ \quad
\Big|I^{N,\eps}(0)-N\overline{I}^\eps(0)\Big| \le 1$;

\smallskip
		
\item there exists two continuous functions $\sbff $, $\ibff$ $: \tor \longrightarrow \R_+$ such that $c\le \sbff(x)\le C $, \; $\ibff(x)\le C$ for all $x\in \tor$, 
$\int_{\tor}\Big[\sbff(x)+\ibff(x)\Big]dx =1$ and $$\sep(0, \xp) = \eps^{-d}\int_{\vx}\sbff(x)dx, \quad \iep(0,\xp) = \eps^{-d}\int_{\vx}\ibff(x)dx\,.$$

\smallskip

		 
\smallskip
	
\item $\{X_j(0)\, , 1\le j\le S^{N,\eps}(0)\}$   and  $\{Y_j(0)\, , 1\le j\le I^{N,\eps}(0)\}$ are two mutually independent collections of i.i.d. random variables satisfying $\dis \P\left(X_j(0)=\xp\right)=\dfrac{\sep(0, \xp)}{\sep(0)}, \; \text{and}  \; \; \P\left(Y_j(0)=\xp\right)=\dfrac{\iep(0, \xp)}{\iep(0)}$ for all $\xp \in \dispace$,  where
$\dis \sep(0) :=\sum_{\xp} \sep(0, \xp)$ and $\dis \iep(0) :=\sum_{\xp} \iep(0, \xp).$
	Moreover $\dis S^{N,\eps}(0,\xp)=\sum_{j=1}^{S^{N,\eps}(0)}\mds_{X_j(0)=\xp}$ and $\dis I^{N,\eps}(0,\xp)=\sum_{j=1}^{I^{N,\eps}(0)}\mds_{Y_j(0)=\xp}$.
\end{enumerate}
\end{assumption}

\begin{assumption}\label{ass2}	
\begin{enumerate}[label=\bf (\roman*)]
\item  We assume that $\dis  \beta_{\eps}^{\xp, \yp}(t)=\beta_{t}(\xp, \vx)$, where $\beta_{t}(x, A)$ is a transition kernel and there exists a constant $\beta^*$ such that $\beta_{t}(x, \tor)\le \beta^*$, for all $x\in \tor$ and $t\ge 0$.
\item there exists  a positive constant $\lambda^* >0$ such that $0\le \lambda_j(t)\le \lambda^*$, for all $j\in\Z\backslash\{0\}$ and $t\ge0$.
\end{enumerate}
\end{assumption}

Under \ref{ass1} and \ref{ass2}, we have the

\begin{theorem}[\textbf{Law  of  Large  Numbers: \boldmath $\mathbf{N}\to \infty$,  $\eps$  being fixed}]\label{llns}~\\
	As $N\rightarrow \infty$, $\displaystyle \left(\overline{S}^{N, \eps}(t,x_\eps),\, \overline{\mathfrak{F}}^{N, \eps}(t,x_\eps), \, \overline{I}^{N, \eps}(t,x_\eps), \; \; t\ge 0,\, x_\eps \in \dispace\right)$ converges in $\mathbf{D}^{3\eps^{-d}}$ in probability, to the unique solution
	$\displaystyle \left(\overline{S}^\eps(t,x_\eps),\, \overline{\mathfrak{F}}^\eps(t,x_\eps), \, \overline{I}^\eps(t,x_\eps), \; \; t\ge 0,\, x_\eps \in \dispace\right)$ of the following system of integral equations
	
	\begin{equation}\label{det}
	\left\{
	\begin{aligned}
	\overline{S}^{\,\eps}(t,x_\eps) &= \overline{S}^{\,\eps}(0,x_\eps)-\int_{0}^{t}\overline{S}^{\,\eps}(s,x_\eps) \overline{\Gamma}^{\,\eps}(s,x_\eps)ds+\int_0^t  \big[\Delta^S_\eps\overline{S}^{\,\eps}\big](s,x_\eps)ds\\
	\overline{\mathfrak{F}}^{\,\eps}(t,x_\eps) &= \overline{\lambda}_{0}(t)\sum_{y_\eps}\overline{I}^{\,\eps}(0,y_\eps)q^{y_\eps,x_\eps}(0,t)+\sum_{y_\eps}^{}\int_{0}^{t}\overline{\lambda}(t-s)\overline{S}^{\,\eps}(s,y_\eps)\overline{\Gamma}^{\,\eps}(s,y_\eps)q^{y_\eps,x_\eps}(s,t)ds\\
	\overline{I}^{\,\eps}(t,x_\eps) &= \overline{I}^{\,\eps}(0,x_\eps)+\int_{0}^{t}\overline{S}^{\,\eps}(s,x_\eps) \overline{\Gamma}^\eps(s,x_\eps)ds+\int_{0}^{t} \big[\Delta^I_\eps\overline{I}^{\,\eps}\big](s,x_\eps)ds,\\
	t\ge 0 ,&  \; \;  \xp \in \dispace, 
	\end{aligned}
	\right.
	\end{equation}
	where
\[\overline{\Gamma}^{\,\eps}(t,x_\eps) = \frac{1}{\big[\,\overline{B}^{\,\eps}(t,x_\eps)\big]^\gamma}\sum_{y_\eps}\beta_\eps^{x_\eps,y_\eps}(t)\overline{\mathfrak{F}}^{\,\eps}(t,y_\eps)\; \; \text{and} \; \; 
	\overline{B}^{\,\eps}(t,x_\eps) = \overline{S}^{\,\eps}(t,x_\eps) +\overline{I}^{\,\eps}(t,x_\eps).\]
\end{theorem}

This Theorem is a special case of Theorem 3.1 in \cite{FPP22}, whose
 proof written for a multi-patch multi-group SIR model is easily adapted to our case.

\section{Limit as $\eps \to 0$ in the deterministic model}

Before letting $\eps$ go  to zero in the limit system  (\ref{det}) extended on the whole space $\tor$, we prove some technical lemmas.

\begin{lemma}\label{bornSI} Let $T>0$. There exists a positive constant $C$ such that
$\big\Vert \overline{S}^\eps(t) \big\Vert_{\infty} \le C$ and $\big\Vert \overline{I}^\eps(t) \big\Vert_{\infty} \le C$, for all $\eps>0$ and $t\in[0\, ,\, T]$.
\end{lemma}

\begin{pr}
Using the Duhamel formula, we have 
$\displaystyle\|\overline{S}^\eps(t)\|_\infty \le \sup_{\xp} \overline{S}^\eps(0,\xp) \le C$.

We now consider the term $\overline{I}^\eps$. 
First using the previous estimate, we obtain 
 \begin{equation*}
 \frac{\overline{S}^\eps(s,x_\eps)}{\left(\overline{B}^\eps(s,x_\eps)\right)^\gamma}=\left(\frac{\overline{S}^\eps(s,x_\eps)}{\overline{B}^\eps(s,x_\eps)}\right)^\gamma \Big[\,\overline{S}^\eps(s,x_\eps)\Big]^{1-\gamma} \le
 C(T,\gamma). 
 \end{equation*}
   Next we have  $\dis \sum_{y_\eps}\beta_\eps^{x_\eps,y_\eps}\overline{\mathfrak{F}}^\eps(s,y_\eps) \le \lambda^*\left\Arrowvert\overline{I}^\eps(s)\right\Arrowvert_\infty\sum_{y_\eps}\beta_\eps^{x_\eps,y_\eps}(s)\le \lambda^*\beta^*\big\Vert \overline{I}^{\eps}(s)\big\Vert_{\infty}.$ 
Thus
\begin{align*}
\left\Arrowvert\overline{I}^\eps(t)\right\Arrowvert_\infty 
&\le \left\Arrowvert\left(T_{I,\eps}(t)\overline{I}^\eps(0)\right)\right\Arrowvert_\infty + \int_0^tT_{I,\eps}(t-s)C\left\Arrowvert\overline{I}^\eps(s)\right\Arrowvert_\infty ds\\
& \le C + C\int_0^t\left\Arrowvert\overline{I}^\eps(s)\right\Arrowvert_\infty ds.
\end{align*}
The second statement then follows from Gronwall's Lemma. \hfill $\square$
\end{pr}

\begin{lemma}\label{bornB} For any $T>0$, there exists $\eps_0$ and $c>0$ such that $\bep(t, \xp)\ge c$, for all $0<\eps\le \eps_0$, $\xp \in \dispace$ and  $0\le t\le T$.
\end{lemma}

\begin{pr} Let $c$ and $C$ be two positive constants such that
 $\displaystyle 0< c \le \frac{\inf_{x_\eps}\sep(0,x_\eps)}{2} \le \frac{C}{2}$, and let $\displaystyle T_c^{\eps} := \inf\{t >0 \, , \; \inf_{x_\eps}\sep(t,x_\eps) < c\}$.
On the interval $[0\,, \,T_c^{\eps}]$, $\sep(t,x_\eps) \ge c$, $\forall x_\eps \in \dispace$. For $t \le T_c^\eps$, we have
\begin{align*}
\overline{\Gamma}^\eps(t,x_\eps) &=\frac{1}{\big[\,\overline{B}^{\,\eps}(t,x_\eps)\big]^\gamma}\sum_{y_\eps}\beta_\eps^{x_\eps,y_\eps}(t)\overline{\mathfrak{F}}^\eps(t,y_\eps)
\le \frac{\lambda^*\beta^*}{c^\gamma}\big\Vert \overline{I}^{\eps}(t)\big\Vert_{\infty} := \overline{c}, 
\\
\overline{S}^\eps(t,x_\eps) &\le \overline{S}^\eps(0,x_\eps)-\overline{c}\int_{0}^{t}\overline{S}^\eps(s,x_\eps)+\int_0^t  \big[\Delta^S_\eps\overline{S}^\eps\big](s,x_\eps)ds.
\end{align*}
Hence $e^{\overline{c}t}\overline{S}^\eps(t,x_\eps)\ge\inf_{y_\eps}\overline{S}^\eps(0,y_\eps)=2c$, and consequently $T^\eps_c\ge\log2/\overline{c}$.
Then for all $0 \le t \le T_c^{\eps}$, we have $e^{\overline{c}t}\overline{S}^\eps(t,x_\eps) \ge 2c$.  So $\dis \overline{S}^\eps(t,x_\eps) \ge 2e^{-\overline{c}t} c  \ge c $ iff  $e^{-\overline{c}t} \ge \frac{1}{2}$

From \ref{ass1} \textbf{(ii)} and the fact that $\overline{\mathbf{I}}(0)\ne 0$,  there exists a ball $B(x_0, \rho)$ and $a>0$ such that $\overline{\mathbf{I}}(y)\ge a$, for all $y\in B(x_0, \rho)$. Let us consider the following ODE $$\dfrac{ d\,u_\eps}{dt}=\nu_I\Delta_{\eps}u_{\eps}, \quad u_{\eps}(0)= a\mds_{B(x_0, \rho)}.$$ We have that $u_{\eps} \longrightarrow u$ in $L^{\infty}\left([0,T]\times\tor\right)$ as $\eps \to 0$, where $u$ is the solution of 
$$\dfrac{d\,u}{dt}=\nu_I\Delta u, \quad u(0)= a\mds_{B(x_0, \rho)}.$$ 
For all $\dfrac{\log 2}{\overline{c}}< t \le T $, there exists a positive constant $\underline{c}$ , such that $u(t,x) \ge 2\underline{c}, \; \forall x\in \tor$. Then, there exists $\eps_0>0$ such that $\forall \eps \le \eps_0$, $\overline{I}^{\eps}(t,\xp) \ge u_{\eps}(t, x_\eps)\ge \underline{c}$ ,\; for all $\dfrac{\log 2}{\overline{c}}< t \le T $.\\
We have shown that $\overline{B}^{\eps}(t,\xp)\ge c\wedge \underline{c}$, for all $0\le t\le T$, $x\in \dispace$, $\eps_0\le \eps$.

\fpr
\end{pr}

We now extend the  solution of the system (\ref{det}) to the whole space $\tor$. So, we define
\begin{align*}
\sbf(t,x)&:=\sum_{\xp} \overline{S}^{\eps}(t, \xp)\mds_{\vx}(x), \ \ibf(t,x):=\sum_{\xp} \overline{I}^{\eps}(t, \xp)\mds_{\vx}(x), \ 
 \fbf(t,x):=\sum_{\xp} \overline{\mathfrak{F}}^{\eps}(t, \xp)\mds_{\vx}(x),\\  \xbf &:=(\sbf\,,\, \fbf\,,\, \ibf).
\end{align*}
\begin{theorem}\label{CDD}
For all $T \ge 0$, $\dis\sup_{0 \le t \le T}\Big\Vert \xbf(t)-\xbff(t)\Big\Vert_{\infty} \longrightarrow 0$ as $\eps\to0$, where $\xbff :=(\sbff\, ,\, \fbff \,,\,\ibff)$ is 
the unique solution of  the following system of parabolic PDE/integral equations.

\begin{equation}\label{detlim}
\left\{
\begin{aligned}
\sbff(t,x) &= \sbff(0,x)-\int_{0}^{t}\sbff(s,x)\overline{\Gamma}(s,x)ds+\int_{0}^{t} \big[\Delta^{\!S} \sbff\,\big](s,x)ds,\\
\fbff(t, x) &= \overline{\lambda}_0(t)\left(T_I(t)\overline{\mathbf{I}}(0)\right)(x)+\int_0^t \overline{\lambda}(t-s)T_I(t-s)\left(\overline{\mathbf{S}}(s)\overline{\Gamma}(s)\right)(x)ds, \\
\ibff(t,x) &=\ibff(0,x)+\int_{0}^{t}\sbff(s,x)\overline{\Gamma}(s,x)ds+\int_{0}^{t} \big[\Delta^{\!I}\ibff\, \big](s,x)ds,\\
\!\!\!\text{\normalfont with}&\;\;\sbff(t,x)\overline{\Gamma}(t,x)= \frac{\sbff(t,x)}{\big[\,\bbff(t,x)\big]^\gamma} \int_{\tor} \fbff(t,y)\beta(x,dy),\
t\ge 0 ,  \; \;  x \in \tor.
\end{aligned}
\right.
\end{equation}
where $T_I$ denotes the semigroup generated by $\nu_I\Delta$.
\end{theorem}

Before proving this theorem, we first establish two Propositions.

\begin{proposition}\label{bornsup} Let $T > 0$. If $(\sbff\, ,\, \fbff \,,\,\ibff)$ is a solution of (\ref{detlim}), then for all $ 0\le t \le T$, there exists $C$, $c > 0$ such that 
	$\big\Vert \sbff(t)\big\Vert_\infty \le C$, $\big\Vert \ibff(t)\big\Vert_\infty \le C$ and $\bbff(t,x) \ge c$ , for all $x\in \tor$.
\end{proposition}
\begin{pr}
The arguments used in the proof of \ref{bornSI} and \ref{bornB} are easy to transpose  to the present situation.
\hfill $\square$
\end{pr}

\begin{remark}\label{rm1}
Let $ \dis \mathscr{H}\left(\sbff, \ibff, \fbff\right)(t,x) := \dfrac{\left[\,\sbff(t,x)\vee 0\right]\wedge C}{\big[\,\bbff(t,x)\vee c\big]^\gamma} \int_{\tor}\beta_t(x,dy)\left[\,\fbff(t,y)\wedge \lambda^* C\right]$ where $C$ is the upper bound in Lemma \ref{bornSI}, and  $c$ the lower bound in Lemma \ref{bornB}.
Note $\left(\sbff \, , \, \ibff \, , \, \fbff\right)$  is a solution  of (\ref{detlim}) iff it is a solution of the following system
\begin{equation}\label{semgr}
\left\{
\begin{aligned}
\sbff(t,x) &= \Big(T_S(t)\sbff(0)\Big)(x)-\int_{0}^{t}\Big(T_S(t-s)\mathscr{H}\left(\sbff(s), \ibff(s), \fbff(s)\right)\Big)(x)ds,\\
\fbff(t, x) &= \overline{\lambda}_0(t)\Big(T_I(t)\ibff(0)\Big)(x)+\int_0^t \overline{\lambda}(t-s)\Big(T_I(t-s)\mathscr{H}\left(\sbff(s), \ibff(s), \fbff(s)\right)\Big)(x)ds, \\
\ibff(t,x) &= \Big(T_I(t)\ibff(0)\Big)(x)+\int_{0}^{t}\Big(T_I(t-s)\mathscr{H}\left(\sbff(s), \ibff(s), \fbff(s)\right)\Big)(x)ds, \ 
0\le t\le T ,  \; \;  x \in \tor .
\end{aligned}
\right.
\end{equation}
Note also that the map $\dis \mathscr{H}: \left(L^{\infty}(\tor)\right)^{3} \longrightarrow L^{\infty}(\tor)$ is bounded and globally Lispchitz.
\end{remark}

\begin{proposition}
The system of equations  (\ref{semgr}) has a unique solution. 
\end{proposition}

\begin{pr} The uniqueness  of the solution uses the contraction character of the semigroups $T_S$ and $T_I$ on $L^{\infty}(\tor)$, and the fact that the map $\mathscr{H}$ is bounded and globally Lispchitz. The existence of the solution can be proved using the Picard iteration procedure. 
\fpr
\end{pr}

We introduce the canonical projection $\mathrm{P}_{\!\!\eps}: L^2(\tor) \longrightarrow \mathrm{H}^{\eps}$ given by
$$
\varphi \longmapsto \mathrm{P}_{\!\!\eps} \varphi(x)=\varepsilon^{-d} \int_{\vx} \varphi(y) d y \quad \text { if } x \in \vx.
$$

\bigskip

\textbf{Proof of Theorem} \ref{CDD}.

Using the fact that the map $\mathscr{H}$ is bounded and globally Lispchitz, we have, provided that $\eps\le\eps_0$,  
\begin{align*}
\Big\Vert \xbf(t)-\xbff(t)\Big\Vert_\infty \le C(\lambda^*, \beta^*) \int_{0}^{t}\Big\Vert \xbf(s)-\xbff(s)\Big\Vert_\infty ds + \pi_\eps(t),
\end{align*}
where $\pi_\eps(t) = \pi^S_\eps(t)+\pi^I_\eps(t)+\pi^{\mathfrak{F}}_\eps(t)$,
with
\begin{align*}
\pi^S_\eps(t) &= \left\Arrowvert T_{S,\eps}(t)\sbf(0) -T_S(t)\sbff(0)\right\Arrowvert_\infty\\
&\hspace{-1cm} + \int_{0}^{t}\left\Arrowvert \mathrm{P}_{\!\!\eps}\left(\frac{\sbff(s)}{\left[\,\bbff(s)\right]^\gamma} \int_{\tor}\fbff(s,y)\beta_s(.,dy)\right)-\frac{\sbff(s)}{\left[\,\bbff(s)\right]^\gamma} \int_{\tor}\fbff(s,y)\beta_s(.,dy)\right\Arrowvert_\infty ds\\
& \hspace{-1cm}+ \int_{0}^{t}\left\Arrowvert T_{S,\eps}(t-s)\mathrm{P}_{\!\!\eps}\left(\frac{\sbff(s)}{\left[\,\bbff(s)\right]^\gamma} \int_{\tor}\fbff(s,y)\beta_s(.,dy)\right)-T_S(t-s)\left(\frac{\sbff(s)}{\left[\,\bbff(s)\right)^\gamma} \int_{\tor}\fbff(s,y)\beta_s(.,dy)\right)\right\Arrowvert_\infty ds,
\end{align*}
$\pi^I_\eps(t)$ is a quantity similar to $\pi^S_\eps(t)$, with $T_{I, \eps}$ (resp. $T_I$, $\ibf$ and $\ibff$) in place of $T_{S, \eps}$ (resp. $T_S$, $\sbf$ and $\sbff$), and 
\begin{align*}
\pi^{\mathfrak{F}}_\eps(t) &= \lambda^*\Big\Vert T_{I,\eps}(t)\ibf(0)-T_I(t)\ibff(0)\Big\Vert_\infty\\
& \hspace{-1cm} + \int_{0}^{t}\Big\Vert \mathrm{P}_{\!\!\eps}\left(\frac{\sbff(s)}{\left[\,\bbff(s)\right]^\gamma} \int_{\tor}\fbff(s,y)\beta_s(.,dy)\right)-\frac{\sbff(s)}{\left[\,\bbff(s)\right]^\gamma} \int_{\tor}\fbff(s,y)\beta_s(.,dy)\Big\Vert_\infty ds\\
& \hspace{-1cm}+ \int_{0}^{t}\Big\Vert T_{I,\eps}(t-s)\mathrm{P}_{\!\!\eps}\left(\frac{\sbff(s)}{\left[\,\bbff(s)\right]^\gamma} \int_{\tor}\fbff(s,y)\beta_s(.,dy)\right)-T_I(t-s)\left(\frac{\sbff(s)}{\left[\,\bbff(s)\right]^\gamma} \int_{\tor}\fbff(s,y)\beta_s(.,dy)\right)\Big\Vert_\infty ds .
\end{align*}
Then from Gronwall’s lemma, 
$\sup_{0 \leq t \leq T}\Big\Vert\overline{\mathbf{X}}^{\,\eps}(t) -\overline{\mathbf{X}}(t)\Big\Vert_\infty \to0$ follows from  $\sup_{0 \leq t \leq T}\pi_\eps(t)\to0$.

Since the maps $x\longmapsto \sbff(0,x)$, $x\longmapsto \ibff(0,x)$ and $x\longmapsto \frac{\sbff(t,x)}{\left[\,\bbff(t,x)\right]^\gamma} \int_{\tor}\fbff(t,y)\beta_t(x,dy)$ are continuous on $\tor$, and the fact that $T_{S,\eps}\longrightarrow T_{S} $ and $T_{I,\eps}\longrightarrow T_{I}$ in $L^{\infty}$ as $\eps \to 0$, then  $\dis\sup_{0 \leq t \leq T}\pi_\eps(t) \longrightarrow 0$,  as $\eps \to 0$  (see Kato \cite{TKTO}, chapter 9, Section 3, Example 3.10).
\fpr

\section{Limit as $N\to \infty$ and $\eps \to 0$}

\medskip

In this section, we extend our stochastic model on the whole space $\tor$ and let both $N\to \infty$ and $\eps \to 0$ in such a way that $N\eps^d\to \infty$. Before stating the main theorem of this section, we first prove some lemmas and propositions.

\begin{lemma}\label{p1}
	There exist two constants $0< c<C$ such that for all $t\ge 0$, $\eps>0$ and $x_\eps \in \dispace$, $$\displaystyle c\eps^d \leq \mathbb{P}(X(t)=x_\eps)\leq C\eps^d .$$
\end{lemma}
\begin{pr} Define 	$\displaystyle u^\eps(t,x_\eps) := \mathbb{P}(X(t)=x_\eps)$. 
	We have that $\displaystyle u^\eps(t,x_\eps) = \left(e^{t\big[\Delta^S_\eps\big]^*}u^\eps_0\right)(x_\eps)$. Using the assumption on the initial condition $\mathbb{P}(X(0)=x_\eps)$, then  $0<c\eps^{d} \leq u^\eps(0, x_\eps) \leq C\eps^{d}$, from which we deduce that   $0<c\eps^{d} \leq e^{t\big[\Delta^S_\eps\big]^*}u^\eps(0, x_\eps) \leq C\eps^{d}$, hence the result.
	\end{pr}

\begin{lemma}\label{lemt0}
	There exits a positive constant $C$ such that for all $0 \le s \le t$, $\eps>0$ and $x_\eps \in \dispace$
	$$\sum_{y_\eps}q^{y_\eps,x_\eps}_{\eps}(s,t) = 1 \quad \text{and} \quad 
\mathbb{P}\left(Y_j(t)=x_\eps\right) \le C\eps^d.$$
\end{lemma}

\begin{pr}
The uniform distribution on $\dispace$ is invariant for the process $Y(t)$. So if we start $Y$ at  time $s$ with the uniform distribution i.e. $\mathbb{P}\left(Y(s)=x_\eps\right) = \eps^{d}$, the law of $Y$ at time $t$ is also the uniform law. But
\begin{align*}
\mathbb{P}\left(Y(t)=x_\eps\right) = \sum_{y_\eps}\mathbb{P}\left(Y(s)=y_\eps\right)q^{y_\eps,x_\eps}_{\eps}(s,t) \; \text{i.e} \; \eps^d =  \eps^d\sum_{y_\eps}q^{y_\eps,x_\eps}_{\eps}(s,t),
\end{align*}
thus
$\; \; \; \dis \sum_{y_\eps}q^{y_\eps,x_\eps}_{\eps}(s,t) = 1.$ Finally
\begin{align*}
\mathbb{P}\left(Y_j(t)=x_\eps\right) 
&= \sum_{y_\eps}\mathbb{P}\left(Y_j(0)=y_\eps\right)q^{y_\eps,x_\eps}_{\eps}(0,t)\\
&\le \sup_{y_\eps}\mathbb{P}\left(Y_j(0)=y_\eps\right)\sum_{y_\eps}q^{y_\eps,x_\eps}_{\eps}(0,t).
\end{align*}
Hence the second result follows from the first one and  Assumption \ref{ass1} $\mathbf{(ii)}$ and $\mathbf{(iii)}$.\hfill $\square$	
\end{pr}

Let define
$\displaystyle\overline{\mathfrak{F}}_0^{N,\eps}(t,x_\eps):=\frac{1}{N} \sum_{j=1}^{I^{N,\eps}(0)}\lambda_{-j}(t)\mds_{Y_j(t)=x_\eps}$ and\;
$\displaystyle\overline{\mathfrak{F}}_0^\eps(t,x_\eps) := \overline{\lambda}_0(t)\sum_{y_\eps}^{}\overline{I}^\eps(0,y_\eps)q^{y_\eps,x_\eps}_{\eps}(0,t)$.

We have the

\begin{lemma}\label{ff0}
	Let us assume that $(N,\eps)\to (\infty, 0)$, in such a way that $N\eps^d\to \infty$.
	Then for all $T>0$,  $$\sup_{0\le t\le T}\E\left(\left\|\overline{\mathfrak{F}}_0^{N,\eps}(t)-\overline{\mathfrak{F}}_0^\eps(t)\right\|_\infty^2\right) \longrightarrow 0,\quad
	\text{as} \; \; (N\, ,\, \eps) \rightarrow (\infty\, ,\, 0).$$	
\end{lemma}
\begin{pr}
 $\displaystyle\overline{\mathfrak{F}}_0^{N,\eps}(t,x_\eps)$  can be  decomposed as follows \\
	$\displaystyle\overline{\mathfrak{F}}_0^{N,\eps}(t,x_\eps)  = \frac{1}{N} \sum_{j=1}^{I^{N,\eps}(0)}\left(\lambda_{-j}(t)-\overline{\lambda}_0(t)\right)\mds_{Y_j(t)=x_\eps}+\overline{\lambda}_0(t)\frac{1}{N}\sum_{j=1}^{I^{N,\eps}(0)}\mds_{Y_j(t)=x_\eps}$.\\
	Let consider the first term. 
	Since $\left(\lambda_{-j}(t)\right)_j$ are independent and identically distributed and independent of $Y_j(t)$, then
	\begin{align*}
	\mathbb{E}\left[\left(\frac{1}{N} \sum_{j=1}^{I^{N,\eps}(0)}\left(\lambda_{-j}(t)-\overline{\lambda}_0(t)\right)\mds_{Y_j(t)=x_\eps}\right)^2\right] 
	&= \frac{1}{N^2} \sum_{j=1}^{I^{N,\eps}(0)} \mathbb{E}\left[\left|\lambda_{-j}(t)-\overline{\lambda}_0(t)\right|^2\mds_{Y_j(t)=x_\eps}\right]\\
	&\le \frac{1}{N^2}C(\lambda^*) I^{N,\eps}(0)\mathbb{P}\left(Y_1(t)=x_\eps\right)
	\le \frac{C(\lambda^*)}{N}.
	\end{align*}
	Now, since
	\begin{align}\label{c1}
	\mathbb{E}\left[\sup_{x_\eps \in \dispace}\left(\frac{1}{N} \sum_{j=1}^{I^{N,\eps}(0)}\left(\lambda_{-j}(t)-\overline{\lambda}_0(t)\right)\mds_{Y_j(t)=x_\eps}\right)^2\right]
	&\le \sum_{x_\eps}\mathbb{E}\left[\left(\frac{1}{N} \sum_{j=1}^{I^{N,\eps}(0)}\left(\lambda_{-j}(t)-\overline{\lambda}_0(t)\right)\mds_{Y_j(t)=x_\eps}\right)^2\right]\nonumber
	\\&\le \frac{C(\lambda^*)}{N}\eps^{-d}\quad  \to0,
	\end{align}
	provided $N\eps^d\to \infty$.
	It remains to show that 
	\[\sup_{x_\eps \in \dispace}\left|\overline{\lambda}_0(t)\frac{1}{N}\sum_{j=1}^{I^{N,\eps}(0)}\mds_{Y_j(t)=x_\eps}-\overline{\lambda}_0(t)\sum_{y_\eps}\overline{I}^\eps(0,y_\eps)q^{y_\eps,x_\eps}_{\eps}(0,t)\right| \longrightarrow 0, \text{ as }(N,\eps) \longrightarrow (\infty,0).\]
	We have
	\begin{align*}
	\frac{1}{N}\sum_{j=1}^{I^{N,\eps}(0)}\mds_{Y_j(t)=x_\eps} = \frac{1}{N}\sum_{j=1}^{I^{N,\eps}(0)}\left[\mds_{Y_j(t)=x_\eps}-\mathbb{P}\left(Y_j(t)=x_\eps\right)\right] + \frac{1}{N}\sum_{j=1}^{I^{N,\eps}(0)}\mathbb{P}\left(Y_j(t)=x_\eps\right).
	\end{align*}
	\begin{align*}
	\mathbb{E}\left\{\left(\frac{1}{N}\sum_{j=1}^{I^{N,\eps}(0)}\left[\mds_{Y_j(t)=x_\eps}-\mathbb{P}\left(Y_j(t)=x_\eps\right)\right]\right)^2 \right\} &= \frac{1}{N^2}\sum_{j=1}^{I^{N,\eps}(0)}\mathbb{E}\left(\Big\vert\mds_{Y_j(t)=x_\eps}-\mathbb{P}\left(Y_j(t)=x_\eps\right)\Big\vert^2 \right)\\
	&\le \frac{C}{N}\,.
	\end{align*}
	It follows that
	\begin{align}\label{c2}
	\mathbb{E}\left\{\sup_{x_\eps \in \dispace}\left(\frac{1}{N}\sum_{j=1}^{I^{N,\eps}(0)}\left[\mds_{Y_j(t)=x_\eps}-\mathbb{P}\left(Y_j(t)=x_\eps\right)\right]\right)^2 \right\} 
	&\le \frac{C}{N\eps^d}\quad \to0,
		\end{align}
	provided $N\eps^d\to 0$.
	
	Since $\overline{\lambda}_0(t)$ is bounded, 
	it  remains to evaluate the quantity  $\displaystyle\frac{1}{N}\sum_{j=1}^{I^{N,\eps}(0)}\mathbb{P}\left(Y_j(t)=x_\eps\right)-\sum_{y_\eps}^{}\overline{I}^\eps(0,y_\eps)q^{y_\eps,x_\eps}_{\eps}(0,t)$.
	We have
	\begin{align*}
	\frac{1}{N}\sum_{j=1}^{I^{N,\eps}(0)}\mathbb{P}\left(Y_j(t)=x_\eps\right)&= \frac{1}{N}\sum_{y_\eps}\sum_{j=1}^{I^{N,\eps}(0)} \mathbb{P}\left(Y_j(0)=y_\eps\right)q^{y_\eps,x_\eps}_{\eps}(0,t), \text{ thus}
	\end{align*}
	\begin{align}\label{c3}
	\sup_{x_\eps}\bigg|\frac{1}{N}\sum_{j=1}^{I^{N,\eps}(0)}\mathbb{P}\left(Y_j(t)=x_\eps\right)-\sum_{y_\eps}^{}\overline{I}^\eps(0,y_\eps)q^{y_\eps,x_\eps}_{\eps}(0,t)\bigg| \nonumber
	&\le \frac{1}{N}\sup_{x_\eps}\sum_{y_\eps}q^{y_\eps,x_\eps}_{\eps}(0,t)\bigg|\sum_{j=1}^{I^{N,\eps}(0)}\mathbb{P}\left(Y_j(0)=y_\eps\right)-N\overline{I}^\eps(0,y_\eps)\bigg|\nonumber \\
	&\le \frac{1}{N}\sup_{x_\eps}\sum_{y_\eps}q^{y_\eps,x_\eps}_{\eps}(0,t)\dfrac{\overline{I}^{\eps}(0,\yp)}{\overline{I}^{\eps}(0)}\bigg|I^{N,\eps}(0)-N\overline{I}^{\eps}(0)\bigg|\nonumber \\
	&\le \dfrac{C}{N}\quad  \longrightarrow 0\,.
	\end{align}
	Combining  (\ref{c1}), (\ref{c2}) and $(\ref{c3})$, we finally have
	\begin{align}
	\sup_{0\le t\le T}\E\left(\sup_{x_\eps  \in \dispace}\left|\frac{1}{N} \sum_{j=1}^{I^{N,\eps}(0)}\lambda_{-j}(t)\mds_{Y_j(t)=x_\eps}-\overline{\lambda}_0(t)\sum_{y_\eps}\overline{I}^\eps(0,y_\eps)q^{y_\eps,x_\eps}_{\eps}(0,t)\right|^2 \right)\longrightarrow 0 \, ,
	\end{align}
	 provided $N\eps^d \rightarrow +\infty$.
	 \fpr\end{pr}
 
 Let $\sigma^{N,\eps}$  be the stopping time defined by 
 \begin{align}\label{ta}
 \sigma^{N,\eps}(\omega) := \inf\left\{t>0 \, ,  \omega \notin A_{t,\delta}\cap B_{t,\delta} \right\},
 \end{align}
 where for all $t \le T$, $\delta>0$, 
 \begin{equation*}
 A_{t,\delta} = \left\{\displaystyle\left\Arrowvert\int_0^t T_{S,\eps}(t-s)d\mathscr{M}_S^{N,\eps}(s)\right\Arrowvert_\infty \le \delta\right\},
 \quad
 B_{t,\delta} = \left\{\left\Arrowvert\int_0^t T_{I,\eps}(t-s)d\widetilde{\mathscr{M}}^{\;N,\eps}_{I}(s)\right\Arrowvert_\infty \le \delta\right\},
 \end{equation*}
  with 
  \begin{align*} 
  \mathscr{M}_S^{N,\eps}(t) &= \sum_{y_\eps\sim x_\eps}^{}\frac{1}{N}M_{S}^{y_\eps,x_\eps}\left(N\int_{0}^{t}\frac{\nu_S}{\eps^2}\overline{S}^{N,\eps}(s,y_\eps)ds\right)-\sum_{y_\eps\sim x_\eps}\frac{1}{N}M_{S}^{x_\eps,y_\eps}\left(N\int_{0}^{t}\frac{\nu_{S}}{\eps^2}\overline{S}^{N,\eps}(s,x_\eps)ds\right),\\
 \dis \wi(t) &=\mathscr{M}_I^{N,\eps}(t)+ \msi(t),\quad\text{ where}
 \end{align*}
\begin{align*}
  \mathscr{M}_I^{N,\eps}(t) &= \sum_{y_\eps\sim x_\eps}\frac{1}{N}M_{I}^{y_\eps,x_\eps}\left(N\int_{0}^{t}\frac{\nu_I}{\eps^2}\overline{I}^{N,\eps}(s,y_\eps)ds\right)-\sum_{y_\eps\sim x_\eps}^{}\frac{1}{N}M_{I}^{x_\eps,y_\eps}\left(N\int_{0}^{t}\frac{\nu_{I}}{\eps^2}\overline{I}^{N,\eps}(s,x_\eps)ds\right),\\
 \dis\msi(t)&=\frac{1}{N} \int_0^t\int_0^\infty\mds_{u\le S^{N,\eps}(s^-,x_\eps)\overline{\Gamma}^{N,\eps}(s^-,x_\eps)} \overline{Q}^{x_\eps}(ds,du).
 \end{align*}

 \medskip
 
  $\overline{Q}^{x_\eps}(ds,du):= Q^{x_\eps}(ds,du)-dsdu$  is the compensated PRM associated with $Q^{x_\eps}_{\eps}(ds,du)$, and we have used the notations
  $$\dis M_S^{\xp, \yp}(t)=P_S^{\xp, \yp}(t)-t,\quad  \dis M_I^{\xp, \yp}(t)=P_I^{\xp, \yp}(t)-t.$$
 
 \bigskip
 
 Let $\dis \bar{\bar{c}}:=\frac{\lambda^*\beta^*\left\Arrowvert\overline{I}^{N,\eps}(t)\right\Arrowvert_\infty}{c^\gamma}$ , where  $c$ stands for the bound in \ref{bornB}.  We define the stopping time $$\tau^{N,\eps} = \inf\left\{t>0 \, , \left\Arrowvert\int_{0}^{t}e^{(t-s)\left(\Delta^S_\eps-\bar{\bar{c}}I_d\right)}d\ws(s)\right\Arrowvert_\infty \ge \frac{c}{8} \right\},$$ where $\dis I_d$ is the identity operator on $\mathrm{H}^{\eps}$, and  $\ws(t,x_\eps):=\mathscr{M}_{S}^{N,\eps}(t,x_\eps)-\msi(t,x_\eps)$.
 
 In the proof of the next Proposition, we shall need the following Lemma.
 \begin{lemma}\label{S0}
 As $ (N , \eps) \to (\infty , 0)$ in such way that $N\eps^d \to \infty$,
 $\big\Vert \overline{S}^{N,\eps}(0, .)-\overline{S}^{\eps}(0, .) \big\Vert_{\infty} \longrightarrow 0 $  in $L^2(\Omega)$ .
  \end{lemma}
  \begin{pr}
 We have
 \begin{eqnarray}
 \overline{S}^{N,\eps}(0,x_\eps)-\overline{S}^{\eps}(0,x_\eps)&=&\dfrac{1}{N}\sum_{j=1}^{S^{N,\eps}(0)}\mds_{X_j=x_\eps} - \P\left(X=x_\eps\right)\overline{S}^{\eps}(0)\n\\
 &=& \overline{S}^{\eps}(0)\dfrac{1}{N \overline{S}^{\eps}(0)}\sum_{j=1}^{S^{N,\eps}(0)}\left[\mds_{X_j=x_\eps} - \P\left(X=x_\eps\right)\right]+\dfrac{\P\left(X=\xp\right)}{N}\left[\overline{S}^{N,\eps}(0)-N\overline{S}^{\eps}(0)\right].\n
 \end{eqnarray}
\begin{eqnarray}
\E\left[\Big\vert\overline{S}^{N,\eps}(0,x_\eps)-\overline{S}^{\eps}(0,x_\eps)\Big\vert^2\right]&\le& \dfrac{2}{N^2}\sum_{j=1}^{S^{N,\eps}(0)} Var\left[\mds_{X=x_\eps}\right]+\dfrac{2\left[\P\left(X=\xp\right)\right]^2}{N^2}\n\\
&\le& \dfrac{\overline{S}^{\eps}(0)}{N}\dfrac{C}{c}\eps^{d} +\dfrac{C\eps^{2d}}{N^2} 
\le  \dfrac{C^\prime}{N} +\dfrac{C\eps^{2d}}{N^2}\n.
\end{eqnarray}
Then
\begin{eqnarray}
\E\left[\sup_{x_\eps  \in \dispace}\big\vert\overline{S}^{N,\eps}(0,x_\eps)-\overline{S}^{\eps}(0,x_\eps)\big\vert^2\right]
&\le & \dfrac{C^\prime}{N\eps^d}+\dfrac{C\eps^{d}}{N^2} \n.
\end{eqnarray}
The result follows.
\hfill $\square$
\end{pr}

 \begin{proposition}\label{ppSI}
 	For all $T > 0$, there exists $C$ such that for $N$ large enough if $t \le \sigma^{N,\eps} \wedge T$, then $\left\Arrowvert\overline{S}^{N,\eps}(t)\right\Arrowvert_\infty \le C$ and $\left\Arrowvert\overline{I}^{N,\eps}(t)\right\Arrowvert_\infty \le C$, for all $\eps>0$. Moreover there exists $\eps_0>0$ and $c_0>0$ such that if $t \le \sigma^{N,\eps}\wedge \tau^{N, \eps}\wedge T$, $ \overline{B}^{N,\eps}(t,x_\eps)\ge c_0$, for all $\xp \in \dispace$, provided $\eps\le \eps_0$.
 \end{proposition}
 \begin{pr}
 Let first treat the term
 $\big\Vert\overline{S}^{N,\eps}(t)\big\Vert_\infty$.\\
Using the Duhamel formula,  we have $$\displaystyle\overline{S}^{N,\eps}(t,x_\eps) \le \left(T_{S,\eps}(t)\overline{S}^{N,\eps}(0,.)\right)(x_\eps)+\int_{0}^{t}\left(T_{S,\eps}(t-s)d\mathscr{M}_S^{N,\eps}(s,.)\right)(x_\eps).$$  Since $\overline{S}^{N,\eps}(0,\xp) \le C$, for all $\xp \in  \dispace$, we obtain that for $t \le \sigma^{N,\eps}\wedge T$,
\[\Vert\overline{S}^{N,\eps}(t)\big\Vert_\infty\le C+\delta\,.\]
 
 We now consider  the term $\left\Arrowvert\overline{I}^{N,\eps}(t)\right\Arrowvert_\infty$. Arguing as in the proof of Lemma \ref{bornSI}, we have for $t \le \sigma^{N,\eps}\wedge T$,
 \begin{align*}
 \left\Arrowvert\overline{I}^{N,\eps}(t)\right\Arrowvert_\infty 
 &\le e^{Ct}\left(C+\sup_{0 \le t \le T}\left\Arrowvert\int_0^tT_{I,\eps}(t-s)d\wi(s)\right\Arrowvert_\infty \right)\\
&\le \left(C+\delta\right)e^{CT}.
 \end{align*}
  
 We finally consider the term $\overline{B}^{N,\eps}(t,x_\eps)$.
 It follows from  Lemma \ref{S0} that  $\big\Vert \overline{S}^{N,\eps}(0, .)-\overline{S}^{\eps}(0, .) \big\Vert_{\infty} \longrightarrow 0 $ and from  Lemma \ref{bornB} that $\overline{S}^{\eps}(0, x_{\eps})\ge c $, for all $\xp \in \dispace$, then for $N$ large enough, $\P\left(\inf_{x_\eps}\overline{S}^{N,\eps}(0, x_{\eps})\ge \frac{c}{2}\right)$ is close to $1$.
 Let $\dis T_c^{N, \eps} = \inf\Big\{t \, ,  \inf_{x_\eps}\overline{S}^{N,\eps}(t,x_\eps) < \dfrac{c}{4}\Big\}$.
 On the interval $[0\, ,T_c^{N,\eps})$, $\overline{S}^{N,\eps}(t,x_\eps) \ge \dfrac{c}{4}$, $\forall x_\eps \in \dispace$. For all $t \le T_c^{N ,\eps}\wedge \sigma^{N,\eps}\wedge T$, we have
 \begin{align*}
 \overline{\Gamma}^{N,\eps}(t,x_\eps) =\frac{1}{\big[\,\overline{B}^{N,\eps}(t,x_\eps)\big]^\gamma}\sum_{y_\eps}\beta_\eps^{x_\eps,y_\eps}(t)\overline{\mathfrak{F}}^{N,\eps}(t,y_\eps)
 \le \frac{4^{\gamma}\lambda^*\beta^*\left\Arrowvert\overline{I}^{N,\eps}(t)\right\Arrowvert_\infty}{c^\gamma} = \bar{\bar{c}}
 \end{align*} and then, if moreover $t\le\tau^{N,\eps}$, 
 \begin{eqnarray}
 \overline{S}^{N,\eps}(t,x_\eps) &\ge& \left(e^{(\Delta_{\eps}^S- \bar{\bar{c}}I_d)t} \overline{S}^{N,\eps}(0)\right)(\xp)+\int_{0}^{t}\left(e^{(t-s)(\Delta_{\eps}^S- \bar{\bar{c}}I_d)}d\ws(s)\right)(\xp)\n\\
 &\ge & \dfrac{c}{2}e^{-\bar{\bar{c}}t} - \dfrac{c}{8}.
 \end{eqnarray}
  We note that $\frac{c}{2}e^{-\bar{\bar{c}}t} \ge \dfrac{c}{4} \quad \text{iff} \quad t \le \frac{\log2}{\bar{\bar{c}}} = T_{\bar{\bar{c}}}$.

   So, on the event $\tau^{N, \eps}\wedge \sigma^{N, \eps} \wedge T\ge T_{\bar{\bar{c}}}$,   $\displaystyle \overline{S}^{N,\eps}(t,x_\eps) \ge \frac{c}{8}$,\; \;  $\forall\, 0\le t \le T_{\bar{\bar{c}}}$\, .
 $$\text{For}\;  t > T_{\bar{\bar{c}}},  \quad \overline{I}^{N,\eps}(t,x_\eps) \ge \left(T_{I,\eps}(t)\overline{I}^{N,\eps}(0)\right)(x_\eps) +\int_{0}^{t}\left(T_{I,\eps}(t-s)d\mathscr{M}_{I}^{N,\eps}(s)\right)(x_\eps).$$ 
 We choose $T>T_{\bar{\bar{c}}} $ arbitrary. We know from the proof of Lemma \ref{bornB} that there exists $\eps_0$ and $\underline{\underline{c}}$ such that $\overline{I}^{\eps}(t,\xp)\ge \underline{\underline{c}}$ for all $\eps \le \eps_0$, $\xp\in \dispace$ and $\dfrac{\log 2}{\bar{\bar{c}}}\le t\le T$. If we now choose $\delta=\dfrac{\underline{\underline{c}}}{2}$ in the definition of $\sigma^{N, \eps}$, we deduce that for any $\eps\le \eps_0$, $\xp \in \dispace$, $T_{\bar{\bar{c}}}\le t\le \sigma^{N, \eps}\wedge T$, \;  $\overline{I}^{N,\eps}(t, \xp)\ge\dfrac{\underline{\underline{c}}}{2}.$
\end{pr} \fpr

From now, we decree that $\sigma^{N, \eps}=0$ whenever $\dis \inf_{x_\eps}\overline{S}^{N,\eps}(0, \xp)<\dfrac{c}{2}$, or $\eps>\eps_0$.

 \begin{lemma}\label{SGama}
 	Given $T > 0$, there exists $C>0$ such that for any $t < \tau^{N,\eps}\wedge\sigma^{N,\eps}$, we have
 	\begin{align}
 	\begin{aligned}
 	\left\Arrowvert\overline{S}^{N,\eps}(t)\overline{\Gamma}^{N,\eps}(t)-\overline{S}^\eps(t)\overline{\Gamma}^\eps(t)\right\Arrowvert_\infty
 	&\le C\Bigg(\left\Arrowvert\overline{S}^{N,\eps}(t)-\overline{S}^\eps(t)\right\Arrowvert_\infty\\
 	&+\left\Arrowvert\overline{\mathfrak{F}}^{N,\eps}(t)-\overline{\mathfrak{F}}^\eps(t)\right\Arrowvert_\infty
 	+ \left\Arrowvert \overline{I}^{N,\eps}(t)-\overline{I}^\eps(t)\right\Arrowvert_\infty\Bigg).
 	\end{aligned}
 	\end{align}
 \end{lemma}

\begin{pr} Note that, using the map $\mathscr{H}$ defined in Remark \ref{rm1}, with a slight modification of the constants, we have
	$$\overline{S}^{N,\eps}(t,x_\eps)\overline{\Gamma}^{N,\eps}(t,x_\eps)-\overline{S}^\eps(t,x_\eps)\overline{\Gamma}^\eps(t,x_\eps)=\mathscr{H}\left(\overline{S}^{N,\eps}, \overline{I}^{N,\eps}, \overline{\mathfrak{F}}^{N,\eps}\right)(t, \xp)-\mathscr{H}\left(\overline{S}^{\eps}, \overline{I}^{\eps}, \overline{\mathfrak{F}}^{\eps}\right)(t, \xp),$$  and the result the follows from the fact that  $\mathscr{H}$  is bounded and globally Lipschitz. 
	\fpr
\end{pr}

We define $\dis \omega^{N,\eps}(t)= \omega_S^{N,\eps}(t)+\omega_I^{N,\eps}(t)+ \omega_{\mathfrak{F}}^{N,\eps}(t)$,
with
\begin{equation}\label{defomega}
\begin{split}
 \omega_S^{N,\eps}(t)&=\left\Arrowvert\overline{S}^{N,\eps}(0)-\overline{S}^\eps(0)\right\Arrowvert_\infty +\left\Arrowvert\int_{0}^{t}T_{S,\eps}(t-s)d\ws(s)\right\Arrowvert_\infty,\\
 \omega_I^{N,\eps}(t)&=\left\Arrowvert\bar{I}^{N,\eps}(0)-\overline{I}^\eps(0)\right\Arrowvert_\infty +\left\Arrowvert\int_{0}^{t}T_{I,\eps}(t-s)d\wi(s)\right\Arrowvert_\infty,\\
\omega^{N,\eps}_\mathfrak{F}(t) &= \left\Arrowvert\overline{\mathfrak{F}}_0^{N,\eps}(t)-\overline{\mathfrak{F}}_0^\eps(t)\right\Arrowvert_\infty + \left\Arrowvert\mathscr{M}_{\mathfrak{F}}^{N,\eps}(t) \right\Arrowvert_\infty,
\end{split}
\end{equation}
where
 \[\mathscr{M}_{\mathfrak{F}}^{\; N,\eps}(t,x_\eps) = \frac{1}{N}\sum_{y_\eps}\int_0^t\int_0^\infty \int_\mathbf{D}\int_\mathbf{D} \lambda(t-s){\mathds{1}}_{u\le S^{N,\eps}(s^-,y_\eps)\overline{\Gamma}^{N,\eps}(s^-,y_\eps)} {\mathds{1}}_{Y^{s,y_\eps}(t)=x_\eps}\overline{Q}^{y_\eps}(ds,du,d\lambda,dY)
 .\]
Note that $\dis \mathscr{M}_{\mathfrak{F}}^{\; N,\eps}$ is not a martingale.

 \begin{lemma}\label{wz}
 	As $(N,\eps) \rightarrow (\infty, 0)$, in such a way that $N\eps^d\to \infty$, 
	\begin{align*}
	\sup_{0\le t\le T}\E\left(\mds_{t \le \sigma^{N,\eps}\wedge\tau^{N,\eps}\wedge T}\,[\omega^{N,\eps}(t)]^2\right)   \to 0.
	\end{align*}
 \end{lemma}

\newcommand{\nep}{\mathtt{H}^{^{\eps}}}

\begin{pr} We shall use the following notation

$$\big\Vert \Phi^{\eps} \big\Vert_{\nep}:=\left[\sum_{x_\eps}\big\vert \Phi^{\eps}_{x_\eps}\big\vert^2\right]^{1/2},$$ for any step function $\Phi^{\eps}$ ($\Phi_{x_\eps}^{\eps}$ denoting the value of $\Phi^{\eps}$  on the cell $\vx$).
	
Thanks to Theorem 2.1 in P. Kotelenez \cite{KOPE}, we have
\begin{align*}
\mathbb{E}\left[ \sup_{t \le \sigma^{N,\eps}\wedge\tau^{N,\eps}\wedge T} \left\Arrowvert\int_{0}^{t}T_{S,\eps}(t-s)d\msi(s)\right\Arrowvert^2_{\nep}\right] 
&\le C \mathbb{E}\left[ \left\Arrowvert\msi(\sigma^{N,\eps}\wedge\tau^{N,\eps}\wedge T)\right\Arrowvert^2_{\nep}\right]\\
&\hspace{-5cm} \le \frac{C}{N}\sum_{x_\eps}\mathbb{E}\left(\int_0^T	 \overline{S}^{N,\eps}(s\wedge\sigma^{N,\eps}\wedge\tau^{N,\eps},x_\eps)\overline{\Gamma}^{N,\eps}(s\wedge\sigma^{N,\eps}\wedge\tau^{N,\eps} ,x_\eps) ds\right).
\end{align*}
Provided $t\le \sigma^{N,\eps}\wedge\tau^{N,\eps}\wedge T$, $\overline{\Gamma}^{N,\eps}(t,x_\eps)\le C(\lambda^*, \beta^*)$ and 
 $\overline{S}^{N,\eps}(t,x_\eps) \le C$. Then 
\begin{align*}
\mathbb{E}\left[ \sup_{t \le \sigma^{N,\eps}\wedge \tau^{N,\eps}\wedge T} \left\Arrowvert\int_{0}^{t}T_{S,\eps}(t-s)d\msi(s)\right\Arrowvert^2_{\nep}\right] \le C(\lambda^*, \beta^*)\frac{1}{N\eps^d}.
\end{align*}
Since the $L^\infty$ norm is bounded by the $\nep$ horm, as $(N, \eps)\to (\infty, 0)$, provided $N\eps^d\to0$,
\begin{align}\label{imt}
\mathbb{E}\left[\sup_{t \le \sigma^{N,\eps}\wedge \tau^{N,\eps}\wedge T}\left\Arrowvert\int_0^t T_{S,\eps}(t-s)d\msi(s)\right\Arrowvert_{\infty}^2\right] \longrightarrow 0 .
\end{align}

The same argument can be used for the term $\displaystyle \left\Arrowvert\int_{0}^{t}T_{S,\eps}(t-s)d\mathscr{M}_{S}^{N,\eps}(s)\right\Arrowvert_\infty$.
We conclude that as $(N,\eps) \longrightarrow (\infty,0)$, in such a way that $N\eps^d\to 0$,  
\begin{equation}\label{convsupS}
\sup_{t \le \sigma^{N,\eps}\wedge \tau^{N,\eps}\wedge T}\omega_S^{N,\eps}(t) \longrightarrow 0\ \text{ in }L^2(\Omega)\,.
\end{equation} 
A similar proof establishes that 
\begin{equation}\label{convsupI}
\sup_{t \le \sigma^{N,\eps}\wedge \tau^{N,\eps}\wedge T}\omega_I^{N,\eps}(t) \longrightarrow 0\ \text{ in }L^2(\Omega)\,.
\end{equation} 
 We now consider $\omega^{N,\eps}_\mathfrak{F}(t)$ 
The convergence  to zero of the first term has been established in Lemma \ref{ff0}. We now consider the second term. We have
\begin{align}\label{inf}
\sup_{t\le T}\mathbb{E}\left(\mathds{1}_{t \le \sigma^{N,\eps}\wedge \tau^{N,\eps} \wedge T }\sup_{x_\eps}\left|\displaystyle\mathscr{M}_{\mathfrak{F}}^{\; N,\eps}(t,x_\eps)\right|^2\right)\n  \\
&\hspace{-5cm}=\frac{1}{N^2}\sup_{t\le T}\mathbb{E}\left[\mathds{1}_{t \le \sigma^{N,\eps}\wedge \tau^{N,\eps}\wedge T}\sup_{x_\eps}\left(\displaystyle\sum_{y_\eps}\int_0^t\int_0^\infty \int_\mathbf{D}\int_\mathbf{D} \lambda(t-s){\mathds{1}}_{u\le S^{N,\eps}(s^-,y_\eps)\overline{\Gamma}^{N,\eps}(s^-,y_\eps)}\right. \right. \n\\ 
& \hspace{-4.5cm}\times \bigg. \bigg. \mds_{Y^{s,y_\eps}(t)=x_\eps}\overline{Q}^{y_\eps}(ds,du,d\lambda,dY)\bigg)^2\bigg]\n\\
&\hspace{-5cm} \le\frac{1}{N^2}\sum_{x_\eps,y_\eps}\mathbb{E}\int_0^{\sigma^{N,\eps}\wedge \tau^{N,\eps}\wedge T} \lambda^2(t-s) S^{N,\eps}(s,y_\eps)\overline{\Gamma}^{N,\eps}(s,y_\eps) q^{y_\eps,x_\eps}_{\eps}(s,t)ds\n\\
&\hspace{-5cm} \le\frac{(\lambda^*)^2}{N}\sum_{x_\eps}\mathbb{E}\left[\int_0^{\sigma^{N,\eps}\wedge \tau^{N,\eps}\wedge T}  \sup_{y_\eps}\left|\overline{S}^{N,\eps}(s,y_\eps)\overline{\Gamma}^{N,\eps}(s, y_\eps)\right|\sum_{y_\eps} q^{y_\eps,x_\eps}_{\eps}(s,t)ds\right]\n\\
&\hspace{-5cm} \le C(\lambda^*)\frac{T}{N\eps^d} .
\end{align}
The result follows. Note that since $\displaystyle\mathscr{M}_{\mathfrak{F}}^{\; N,\eps}(t,x_\eps)$ is not a martingale, the result for  $\omega^{N,\eps}_\mathfrak{F}(t)$ is weaker than \eqref{convsupS} and \eqref{convsupI}.
\hfill $\square$
\end{pr}
 Lemma \ref{wz} clearly implies
\begin{lemma}\label{intwz}
	As $(N,\eps) \longrightarrow (\infty,0)$ in such  way that $N\eps^d \rightarrow \infty$, $\displaystyle\mathds{1}_{t \le \sigma^{N,\eps}\wedge\tau^{N,\eps}\wedge T}\int_0^t\omega^{N,\eps}(s)ds   \longrightarrow 0$ in probability.
\end{lemma}
It remains to establish the next result.
\begin{lemma}\label{sigmaPT}
	As $(N,\eps) \rightarrow (\infty,0)$ ,  $\mathbb{P}\left(\sigma^{N,\eps}<T\right) \longrightarrow 0$ and  $\mathbb{P}\left(\tau^{N,\eps}<T\right) \longrightarrow 0$ .
\end{lemma}

\begin{pr}
	We have
	\begin{equation}\label{sigmaT}
	\begin{aligned}
	\mathbb{P}\left(\sigma^{N,\eps}<T\right) &\le \mathbb{P}\left(\sup_{t \le \sigma^{N,\eps}\wedge T}\left\Arrowvert\int_0^t T_{S,\eps}(t-s)d\mathscr{M}_S^{N,\eps}(s)\right\Arrowvert_\infty \ge \delta/2\right)\\
	&\quad + \mathbb{P}\left(\sup_{t \le\sigma^{N,\eps}\wedge  T}\left\Arrowvert\int_0^t T_{I,\eps}(t-s)d\wi(s)\right\Arrowvert_\infty \ge \delta/2\right).
	\end{aligned}
	\end{equation}
	
	We consider the second term only. The first one is treated similarly.
	 $$\left\Arrowvert\int_0^t T_{I,\eps}(t-s)d\wi(s)\right\Arrowvert_\infty \le  \left\Arrowvert\int_{0}^{t}T_{I,\eps}(t-s)d\msi(s)\right\Arrowvert_\infty+\left\Arrowvert\int_{0}^{t}T_{I,\eps}(t-s)d\mathscr{M}_{I}^{N,\eps}(s)\right\Arrowvert_\infty,$$
	from Proposition 3.2 of \cite{MTE20}, we have
	\begin{align}\label{estimI}
	\mathbb{P}\left(\sup_{t \le \sigma^{N,\eps}\wedge T}\left\Arrowvert\int_0^t T_{I,\eps}(t-s)d\mathscr{M}_{I}^{N,\eps}(s)\right\Arrowvert_\infty \ge \frac{\delta}{2}\right) &\le 4\eps^{-d-2}\exp\left(-\mathtt{a}\frac{\delta^2}{16}N\right)
	\end{align}
	Since we assume that $N\eps^d \longrightarrow 0$,  the right side, hence also the left hand side of \eqref{estimI} tends to $0$.	By Chebyshev's inequality, we have
	\begin{align*}
	\mathbb{P}\left(\sup_{t \le \sigma^{N,\eps}\wedge T}\left\Arrowvert\int_0^tT_{I,\eps}(t-s)d\msi(s)\right\Arrowvert_{\nep} \ge \frac{\delta}{2}\right) &\le \frac{4}{\delta^2}\mathbb{E}\left[\sup_{t \le \sigma^{N,\eps}\wedge T}\left\Arrowvert\int_0^t T_{I,\eps}(t-s)d\msi(s)\right\Arrowvert_{\nep}^2\right].
	\end{align*}
	The right hand side tends to $0$ as shown in the proof of Lemma \ref{wz}. Since the $L^\infty$ norm is bounded by the $\nep$ norm, this finishes the proof that $\mathbb{P}\left(\sigma^{N,\eps}<T\right)\to0$. A similar proof establishes the same result for $\tau^{N,\eps}$.
	\fpr
\end{pr}

We now extend our stochastic process to the whole space $\tor$. So, we define

\[\ssbf(t,x):=\sum_{\xp} \overline{S}^{\eps}(t, \xp)\mds_{\vx}(x), \quad \sibf(t,x):=\sum_{\xp} \overline{I}^{\eps}(t, \xp)\mds_{\vx}(x) \] 
\[\sbbf(t,x):=\sum_{\xp} \bbff^{\eps}(t, \xp)\mds_{\vx}(x), \quad \sfbf(t,x):=\sum_{\xp} \overline{\mathfrak{F}}^{\eps}(t, \xp)\mds_{\vx}(x) \] 
and set $ \sxbf :=(\ssbf\,,\, \sfbf\,,\, \sibf)$.

\bigskip

%

 Let us recall the following Gronwall's lemma.
 \begin{lemma}\label{lgw}
 	Let $\phi$ and $\psi$ be two nonegative Borel measurable locally bounded functions on an interval $[0 ,T)$, with $T< \infty$ and $C$ a non-negative constant. If for all $t \in [0,T)$, the following inequality is satisfied :
 	\begin{align}\label{gw2}
 	\phi(t) \le C \int_0^t\phi(s)ds+\psi(t),
 	\end{align} 
 	then
 	$\displaystyle\phi(t) \le C\int_0^t e^{C(t-s)}\psi(s)ds+\psi(t)$   for all  $t \le T$.
 \end{lemma}

\begin{theorem}\label{LLN}
	 Let us assume that $(N,\eps)\to (\infty, 0)$, in such a way that $N\eps^d\to \infty$. Then $(N,\eps)\to (\infty, 0)$
	\begin{align}
	\left\Arrowvert\sxbf(t)-\xbf(t)\right\Arrowvert_\infty \longrightarrow 0, \; \text{in probability},\; \; \forall \, t \ge 0.
	\end{align}	
\end{theorem}

 \begin{pr} Since $\dis \left\Arrowvert\sxbf(t)-\xbf(t)\right\Arrowvert_\infty=\left\Arrowvert\xepn(t)-\xep(t)\right\Arrowvert_\infty$, it suffices to show that  
 $$\left\Arrowvert\xepn(t)-\xep(t)\right\Arrowvert_\infty  \longrightarrow 0, \; \; \text{in probability, for all $t\ge 0$}.$$ 
 We first consider
 \begin{align*}
 \overline{\mathfrak{F}}^{N,\eps}(t,x_\eps)&=\frac{1}{N} \sum_{j=1}^{I^{N,\eps}(0)}\lambda_{-j}(t)\mds_{Y_j(t)=x_\eps}+\sum_{y_\eps}\int_0^t \overline{\lambda}(t-s)\sepn(s,y_\eps)\overline{\Gamma}^{N,\eps}(s,y_\eps) q^{y_\eps,x_\eps}_{\eps}(s,t)ds+\mathscr{M}_{\mathfrak{F}}^{\; N,\eps}(t,x_\eps),\\
 \overline{\mathfrak{F}}^\eps(t,x_\eps) &= \overline{\lambda}_0(t)\sum_{y_\eps}\overline{I}^\eps(0,y_\eps)q^{y_\eps,x_\eps}_{\eps}(0,t)+\sum_{y_\eps}^{}\int_{0}^{t}\overline{\lambda}(t-s)\overline{S}^\eps(s,y_\eps)\overline{\Gamma}^\eps(s,y_\eps)q^{y_\eps,x_\eps}_{\eps}(s,t)ds.
 \end{align*}
 
 Exploiting Lemma  \ref{SGama}, we have the following: for all  $t \le \sigma^{N,\eps}\wedge\tau^{N, \eps}$
 \begin{align}\label{zzz1}
 \left\|\fepn(t)-\fep(t)\right\|_{\infty} &\le \omega_{\mathfrak{F}}^{N, \eps}(t)  + C\int_0^t\Bigg(\left\Arrowvert\overline{S}^{N,\eps}(s)-\overline{S}^\eps(s)\right\Arrowvert_\infty+\left\Arrowvert\overline{\mathfrak{F}}^{N,\eps}(s)-\overline{\mathfrak{F}}^\eps(s)\right\Arrowvert_\infty\n\\
 &+\left\Arrowvert \overline{I}^{N,\eps}(s)-\overline{I}^\eps(s)\right\Arrowvert_\infty\Bigg)ds.
 \end{align} 
By writing $\; \overline{S}^{N,\eps}(t,x_\eps)-\overline{S}^\eps(t,x_\eps)$  and $\overline{I}^{N,\eps}(t,x_\eps)-\overline{I}^\eps(t,x_\eps)$ in their mild semigroup form, 
 and using estimates in Lemmas \ref{bornSI}, \ref{bornB}, \ref{ppSI} and \ref{SGama}, we obtain, for $t\le \sigma^{N,\eps}\wedge\tau^{N,\eps}\wedge T$ 
 \begin{align}
 \begin{aligned}
 \left\Arrowvert\xepn(t)-\xep(t)\right\Arrowvert_\infty& \le C\int_{0}^{t}\left\Arrowvert\xepn(s)-\xep(s)\right\Arrowvert_\infty ds + \omega^{N,\eps}(t).
 \end{aligned}
 \end{align}
 
 Then, it follows from Gronwall's Lemma   (\ref{lgw}) that
 \begin{equation}\label{zw}
 \begin{aligned}
 \left\Arrowvert\xepn(t)-\xep(t)\right\Arrowvert_\infty
 &\le C\int_{0}^{t}e^{C(t-s)}\omega^{N,\eps}(s) ds+\omega^{N,\eps}(t)\\
 &\le Ce^{Ct}\int_{0}^{t}\omega^{N,\eps}(s)ds + \omega^{N,\eps}(t), \quad \forall\;  t \le \sigma^{N,\eps}\wedge\tau^{N,\eps} .
 \end{aligned}
 \end{equation}
 Consequently using lemmas \ref{wz}, \ref{intwz} and \ref{sigmaPT}, for any $t > 0$,  as $(N,\eps) \rightarrow (\infty,0)$, in such a way that $N\eps^d \longrightarrow \infty$,
 \begin{align*}
 \left\Arrowvert\xepn(t)-\xep(t)\right\Arrowvert_\infty \longrightarrow 0 \; \; \text{ in probability}, \;  \forall t\ge 0.
 \end{align*}
\hfill $\square$
 \end{pr}

 \begin{theorem}\label{LLNSUP}
 	For all $T>0$, as $(N,\eps) \longrightarrow (\infty,0)$ in such  a way that $N\eps^d \rightarrow \infty$, we have,
 	\begin{align*}
 	\sup_{0\le t \le T}\Bigg(\left\Arrowvert \ssbf(t)-\sbf(t)\right\Arrowvert_\infty+ \left\Arrowvert\sibf(t)-\ibf(t) \right\Arrowvert_\infty\Bigg) \longrightarrow 0 ~~\text{in  probability}.
 	\end{align*}
 \end{theorem}
\begin{pr} 
In the proof of the theorem \ref{LLN}, we have established the following:	
\begin{align}
\begin{aligned}
\left\Arrowvert\overline{S}^{\, N,\eps}(t)-\overline{S}^{\,\eps}(t)\right\Arrowvert_\infty \le C\int_{0}^{t}\left\Arrowvert\overline{X}^{\,N,\eps}(s)-\overline{X}^{\,\eps}(s)\right\Arrowvert_\infty ds +\omega_S^{N,\eps}(t)\\
\left\Arrowvert\overline{I}^{N,\eps}(t)-\overline{I}^\eps(t)\right\Arrowvert_\infty \le C\int_{0}^{t}\left\Arrowvert\overline{X}^{\,N,\eps}(s)-\overline{X}^{\,\eps}(s)\right\Arrowvert_\infty ds + \omega_I^{N,\eps}(t) .
\end{aligned}
\end{align}
	It follows that
\begin{align*}
\sup_{0\le t \le \sigma^{N,\eps}\wedge \tau^{N,\eps}\wedge T}\left\Arrowvert\overline{S}^{N,\eps}(t)-\overline{S}^\eps(t)\right\Arrowvert_\infty &\le \sup_{0\le t \le \sigma^{N,\eps}\wedge \tau^{N,\eps}\wedge T}C\int_{0}^{t}\left\Arrowvert\overline{X}^{N,\eps}(s)-\overline{X}^\eps(s)\right\Arrowvert_\infty ds \\
&+\sup_{0\le t \le \sigma^{N,\eps}\wedge \tau^{N,\eps}\wedge T}\omega_S^{N,\eps}(t).
\end{align*}
On the other hand, from \eqref{zw}, for all $t\le \sigma^{N,\eps}\wedge\tau^{N,\eps}$, 
\begin{equation}
\begin{aligned}
\left\Arrowvert\overline{X}^{\, N,\eps}(t)-\overline{X}^{\,\eps}(t)\right\Arrowvert_\infty \le Ce^{Ct}\int_{0}^{t}\omega^{N,\eps}(s)ds + \omega^{N,\eps}(t).
\end{aligned}
\end{equation}
So we deduce from  Lemmas \ref{wz},  \ref{intwz} and \ref{sigmaPT} and \eqref{convsupS} that $$\sup_{0\le t \le T}\left\Arrowvert\overline{S}^{N,\eps}(t)-\overline{S}^\eps(t)\right\Arrowvert_\infty \longrightarrow 0 \; \; \text{in probability as} \; (N,\eps) \longrightarrow (\infty,0),$$ 
and the same is true for $\overline{I}^{N,\eps}(t)-\overline{I}^\eps(t)$.
Thus the claim follows.
\fpr
\end{pr}

We can now state our main result.

\begin{theorem}
	For all $T>0$, as $(N,\eps) \longrightarrow (\infty,0)$ in such  a way that $N\eps^d \rightarrow \infty$, we have,
$$\forall\, t\in [0, T], \quad \Big\Vert\overline{\mathbf{F}}^{N, \eps}(t)- \overline{\mathbf{F}}(t)\Big\Vert_{\infty}\longrightarrow 0, \quad \text{in probability},$$
and 
\begin{align*}
\sup_{0\le t \le T}\Bigg(\left\Arrowvert \ssbf(t)-\overline{\mathbf{S}}(t)\right\Arrowvert_\infty+ \left\Arrowvert\sibf(t)-\overline{\mathbf{I}}(t) \right\Arrowvert_\infty\Bigg) \longrightarrow 0 ~~\text{in  probability}
\end{align*} as $(N, \eps)\to (\infty, 0)$ in such a way that $N\eps^d \to \infty$.
\end{theorem}

\begin{pr}
By using the triangle inequality, the first statement follows from Theorem \ref{CDD} and Theorem \ref{LLN}, 
and the second statement  from Theorem \ref{CDD} and Theorem \ref{LLNSUP}. \fpr
\end{pr}

%

\end{document}